\documentclass[11pt]{amsart}
\usepackage{amssymb,latexsym,amsmath,amsthm,mathrsfs,geometry,color}

\usepackage{fullpage}
\numberwithin{equation}{section}
\usepackage[shortlabels]{enumitem}
\usepackage{hyperref}
\hypersetup{
    colorlinks=true,
    linkcolor=black,
    filecolor=magenta,
    urlcolor=cyan,
    citecolor=black
    }

\newtheorem{theorem}{Theorem}
\newtheorem{lem}{Lemma}

\newtheorem*{theorem*}{Theorem}
\newtheorem*{conjecture*}{Conjecture}
\theoremstyle{remark}
\newtheorem{remark}{Remark}
\theoremstyle{plain}

\newtheorem*{densityhyp*}{Zero-Density Hypothesis (ZDH)}%
\newtheorem*{Sdensityhyp*}{Strong Zero-Density Hypothesis (Strong ZDH)}%
\newtheorem*{prop*}{Property $\mathcal{B}_{b}(m)$}
\newtheorem*{MT}{Montgomery Theorem (MT)}

\begin{document}
\title[Pair Correlation]{Pair Correlation of Zeros of the Riemann Zeta Function I: Proportions of Simple Zeros and Critical Zeros}

\author[Baluyot]{Siegfred Alan C. Baluyot}
\address{Mathematics Department\\ East Carolina University \\ Greenville, NC 27858}
\email{baluyots24@ecu.edu}

\author[Goldston]{Daniel Alan Goldston}
\address{Department of Mathematics and Statistics, San Jose State University, San Jos\'e, CA 95192}
\email{daniel.goldston@sjsu.edu}

\author[Suriajaya]{Ade Irma Suriajaya}
\address{Department of Mathematics and Statistics, San Jose State University, San Jos\'e, CA 95192
\newline \indent Faculty of Mathematics, Kyushu University, Fukuoka 819-0395, Japan}
\email{adeirmasuriajaya@math.kyushu-u.ac.jp}

\author[Turnage-Butterbaugh]{Caroline L. Turnage-Butterbaugh}
\address{Carleton College, Northfield, MN 55057}
\email{cturnageb@carleton.edu}


\dedicatory{Dedicated to Hugh Montgomery, whose insight continues to enlighten and surprise.}

\begin{abstract}
Assuming the Riemann Hypothesis (RH), Montgomery proved a theorem in 1973 concerning the pair correlation of zeros of the Riemann zeta-function and applied this to prove that at least $2/3$ of the zeros are simple. In this paper, we investigate the versatility of the pair correlation method and show, for the first time, that it can be used to prove results on the \textit{horizontal distribution} of zeros of the Riemann zeta-function. 

In earlier work we showed how to remove RH from Montgomery's theorem and, in turn, obtain results on simple zeros assuming conditions on the zeros that are weaker than RH. Here we assume a more general condition, namely that all the zeros $\rho = \beta +i\gamma$ with $T<\gamma\le 2T$ are in a narrow vertical box centered on the critical line with width $b/\log T$. We prove that under this assumption with $b=0.3185$ that at least $2/3$ of zeros are simple and on the critical line.
\end{abstract}
\date{June 8, 2026}

\maketitle
\section{Introduction and Statement of the Main Results}

Let $s=\sigma + it$ with $\sigma, t \in \mathbb{R}$. A central problem in analytic number theory is to quantify the proportion of nontrivial zeros of the Riemann zeta-function $\zeta(s)$ that are simple and on the critical line $\Re(s)=1/2$.
In this paper, we only discuss nontrivial zeros of $\zeta(s)$
and we denote them by $\rho = \beta +i\gamma$, where $\beta, \gamma \in \mathbb{R}$ with $\gamma\neq0$, and the multiplicity of this zero is denoted $m_\rho$.
Throughout we refer to these zeros simply as \lq\lq zeros.\rq\rq

In his now famous paper on the pair correlation of the zeros of the Riemann zeta-function, Montgomery \cite{Montgomery73} introduced a method, conditional on RH, designed to study the vertical distribution of the zeros of $\zeta(s)$. Under RH, Montgomery \cite{Montgomery73} showed asymptotically at least $2/3$ of the zeros are simple, and soon after, Montgomery and Taylor \cite{Montgomery74} improved the proportion to $67.25\%$. In 2020, Chirre, Gon\c{c}alves, and de Laat \cite{CGL2020} used Montgomery's method to show that, under RH, more than $67.92\%$ of the zeros are simple. 

In 1974, Levinson \cite{Levinson74} introduced an unconditional method designed to study critical zeros, and he proved that more than $1/3$ of the zeros are on the critical line, which he improved soon thereafter \cite{Levinson75} to more than $0.3474$. Heath-Brown \cite{Heath-Brown79} and Selberg (unpublished) subsequently showed that more than $34.74\%$ of the zeros are on the critical line \textit{and} simple. As in the case of Montgomery's method, Levinson's method is still actively used. Indeed, in 2020 Pratt, Robles, Zaharescu, and Zeindler \cite{Pra20} proved that more than $41.72\%$ of the zeros are on the critical line and that more than $40.75\%$ of the zeros are on the critical line and simple. For a detailed discussion of the rich history of Levinson's method and its modifications, we refer the reader to the introduction of \cite{Pra20}.

In 1998, Conrey, Ghosh, and Gonek introduced a method which uses the discrete, mollified moments of $\zeta'(\rho)$ to study simple zeros. They showed on RH and an additional hypothesis that at least $19/27=.\overline{703}$ of the zeros are simple. Later Bui and Heath-Brown \cite{BuiHB} showed that this result holds on RH alone.

\subsection{Montgomery's method and critical zeros} In the recent paper \cite{BGST-PC}, we returned to Montgomery's pair correlation method and under an assumption weaker than the Riemann Hypothesis, showed that at least 61.7\% of the zeros are simple. There, we wrote \lq\lq The method of proof neither requires nor provides any information on whether any of these zeros are or are not on the critical line where $\beta=1/2$." In this paper, however, we show that this statement is incorrect. In particular, assuming that the zeros of $\zeta(s)$ are contained in a narrow vertical box centered on the critical line, we prove the new result that the pair correlation method yields at least $67.25\%$ of the zeros are both simple and on the critical line. 

\begin{remark}
Comparing numerical results, we see that the pair correlation method obtains a stronger proportion of critical zeros than does Levinson's method, with the obvious caveat that the pair correlation method requires an assumption on the zeros of $\zeta(s)$ and Levinson's method is unconditional. On the other hand, the arguments presented here reveal a novel observation: the pair correlation method, which was developed to study \textit{vertical distribution} of zeros of $\zeta(s)$ can be used to prove results on the \textit{horizontal distribution} of zeros of $\zeta(s)$ as well. 
\end{remark}

\subsection{Statements of results} For $b>0$, let $B_b$ be the thin vertical box centered on the critical line defined by
\begin{equation} \label{B_0}
B_b := \left\{\,s=\sigma+it\,\mathrel{\Big|}\, \frac12 -\frac{b}{2\log T}< \sigma < \frac12 + \frac{b}{2\log T}, \ T<t\le 2T\right\}, \end{equation} 
and let $N(B_b)$ denote the number of zeros of $\zeta(s)$ in $B_b$. Let $N_s(B_b)$ count the number of zeros in $B_b$ that are simple, and let $N_0(B_b)$ count the number of zeros in $B_b$ that are on the critical line. Finally, let $N_0^s(B_b)$ denote the number of zeros in $B_b$ that are simple and on the critical line. Note that in each of these settings, we count zeros with multiplicity.

\begin{theorem}\label{thm1} Assume that, for all sufficiently large $T$, all the zeros $\rho = \beta +i\gamma$ of $\zeta(s)$ with $T<\gamma \le 2T$ are in $B_b$. Then we have, if $b\to 0$ as $T\to \infty$,

\begin{equation}\label{thm1eq1} N^s_0(B_b) = \sum_{\substack{\rho \\ T<\gamma\le 2T\\ \beta =\frac12\ \text{and} \ \rho \text{ simple}}} 1 \ge \left(\frac23 + o(1)\right) \sum_{\substack{\rho \\ T<\gamma\le 2T}}1 = \left(\frac23 + o(1)\right)N(B_b), \end{equation}
where the summation over the zeros $\rho = \beta +i\gamma $ of $\zeta(s)$ counts with multiplicity.
\end{theorem}
\begin{remark} Using Theorem \ref{thm3} below, we also have \begin{equation}\label{thm1eq2} \max(N_0(B_b), N_s(B_b)) \ge \frac12(N_0(B_b)+ N_s(B_b)) \ge \left(\frac56 + o(1)\right)N(B_b). \end{equation}
Note that the inequality \eqref{thm1eq1} implies $\min(N_0(B_b), N_s(B_b)) \ge N^s_0(B_b)\ge (\frac23 + o(1))N(B_b)$, so that \eqref{thm1eq2} shows that at least one of $N_0(B_b)$ or $N_s(B_b)$ is $\ge (\frac56 + o(1))N(B_b)$.
Furthermore, if $S$ is the set of simple zeros, $C$ the set of critical zeros, and $N_{S\cup C}(B_b)$ the number of zeros that are simple or critical or both in $B_b$, Soundararajan pointed out to us that our method allows us to get
\[N_{S\cup C}(B_b) \ge \left(\frac89 + o(1)\right)N(B_b).\]
\end{remark}

The results in Theorem \ref{thm1} also hold for small fixed values of $b$. This requires the use of a special complex kernel we call a Tsang kernel, together with some easy numerical calculations using Mathematica. In stating our results here, we utilize the asymptotic formula for the number of zeros $\rho$ with $T<\gamma\le 2T$ given in \eqref{zerosumnotation}. 

\begin{theorem}\label{thm2} Assume that, for all sufficiently large $T$, all the zeros $\rho = \beta +i\gamma$ of $\zeta(s)$ with $T<\gamma \le 2T$ are in $B_b$. Then as $T\to \infty$, we have
\[
N^s_0(B_{0.3185}) \ge \Big( 0.66666908 + o(1)\Big)\, \frac{T}{2\pi}\log T,
\]
which recovers the $2/3$ proportion obtained in Theorem \ref{thm1}, and
\[
N^s_0(B_{0.001})\ge \Big(0.67250064 + o(1)\Big)\, \frac{T}{2\pi}\log T,
\]
which recovers Montgomery and Taylor's $67.25\%$ proportion.
\end{theorem}

In the last section of the paper we include Table \ref{table:2} which gives additional results obtained for various values of $b$ using two different Tsang kernels. The table exhibits the fact that the results deteriorate as $b$ increases, and the method ultimately fails when $b\ge 4.2$. 

We use the simpler Tsang  kernel in the proof of Theorem \ref{thm1} here, however in \cite{GS26} it is shown that an ordinary Fej{\'e}r kernel extended to complex values can be used to prove Theorem \ref{thm1} when $b\to 0$ as $T\to \infty$.

\begin{remark} Concerning our assumption on $B_b$ in Theorems \ref{thm1} and \ref{thm2}, we can replace the requirement there are no zeros outside $B_b$ for $T<\gamma \le 2T$ with a suitably strong Zero Density Hypothesis, as was done in \cite[Theorem 3]{BGST-PC}. Also in this region of the critical strip, Brian Conrey\footnote[2]{Zeros Near the Critical Line, unpublished 2024 note.} has proved unconditionally the striking result that more than $99.9\%$ of the zeros with $T<\gamma \le 2T$ are in $B_b$ with $b=16$, or equivalently within a horizontal distance $8/\log T$ of the critical line. 
\end{remark}
\subsection{Horizontal Multiplicity} Soundararajan suggested to us a simplification in our original proofs by introducing a quantity we call the \lq\lq Horizontal Multiplicity." For $g$ a fixed real number, the horizontal multiplicity $H(g)$ is the number of zeros on the horizontal line $s=\sigma +ig$, $0<\sigma <1$, counted with multiplicity.
Assuming RH, then $H(g)$ is either zero or the multiplicity of the zero $\rho = 1/2+ig$. If there are zeros off the critical line, then they come in pairs symmetric to the critical line, and $H(g)$ is the sum of the multiplicities of all the zeros on the horizontal line. Moreover, since $H(g)$ does not depend on the position of the zeros on the horizontal line, $H(g)$ returns the same value whether there are multiple zeros on and off the critical line or a just single zero on the critical line (of multiplicity $H(g)$).

We now let 
\begin{equation} \label{HorMult} N^{\circledast}(T) := \sum_{\substack{\rho,\rho' \\ T<\gamma,\gamma' \le 2T \\ \gamma=\gamma'}}1 = \sum_{\substack{\rho \\ T<\gamma \le 2T }}H(\gamma).\end{equation}
Here we count zeros with multiplicity. The first sum is a double sum over $\rho$ and $\rho'$ and is the mechanism by which the pair correlation method detects multiple zeros. Using an elementary argument, we prove in the next section how a strong enough upper bound on $N^{\circledast}(T)$ gives lower bounds for the number of both simple and critical zeros.
\begin{theorem}\label{thm3} Suppose there exists a constant $\mathbf{C}\ge 1$ and that, as $T\to \infty$, 
\begin{equation} \label{thm3input}
N^{\circledast}(T)
\le \Big(\mathbf{C}+o(1)\Big) \frac{T}{2\pi}\log{T}. \end{equation}
Then asymptotically,
\begin{enumerate}[label={(\roman*)},itemsep=4pt]
\item \label{thm3-i} the proportion of zeros of $\zeta(s)$ which are simple and on the critical line is $\ge 2-\mathbf{C}$,
\item \label{thm3-ii} the average of the proportions of simple zeros and of critical zeros of $\zeta(s)$ is $\ge (3-\mathbf{C})/2$,
\item \label{thm3-iii} the proportion of zeros of $\zeta(s)$ which are either simple or critical or both is $\ge (4-\mathbf{C})/3$.
\end{enumerate}
In order to obtain positive proportions of the zeros, note that \ref{thm3-i} requires $1\le\mathbf{C}<2$, \ref{thm3-ii} requires $1\le\mathbf{C}<3$, and \ref{thm3-iii} requires $1\le\mathbf{C}<4$.
\end{theorem}

\subsection{Outline for the paper}

In Section \ref{sec2} we prove the elementary Theorem \ref{thm3} and further discuss horizontal multiplicity. In Section \ref{sec3} we state and prove an unconditional form of Montgomery's Theorem for his $F$ function. In Section \ref{sec4} we define the Tsang kernel and obtain a formula for evaluating the double sum over zeros of the Tsang kernel evaluated at differences of zeros. In Section \ref{sec5} we prove the main properties of the Tsang kernel. Finally, in Section \ref{sec6} we obtain the constant needed in Theorem \ref{thm3} to bound the average horizontal multiplicity, and then Theorems \ref{thm1} and \ref{thm2} follow immediately from Theorem \ref{thm3}.


\section{Horizontal Multiplicity and Proof of Theorem \ref{thm3}}
\label{sec2}

\subsection{Counting zeros with multiplicity vs weighting zeros by multiplicity}

Let $m_\rho$ denote the multiplicity of the zero $\rho$ of $\zeta(s)$. Consider the set of zeros of $\zeta(s)$ given by
\[Z :=\{ \rho \,:\, \zeta(\rho)=0 \}. \]
Since $Z$ is a set, each zero is counted only once, and thus both simple zeros and multiple zeros are each counted once. We say $Z$ is the set of \textit{distinct} zeros.

On the other hand, zeros of analytic functions are counted by the argument principle with their multiplicity. Thus, it is natural to consider the multiset $\mathcal{Z}$ of zeros, where $\rho \in \mathcal{Z}$ if and only if $\zeta(\rho)=0$, and $\mathcal{Z}$ contains $m_\rho$ copies of $\rho$. With this notation, we define $N(T)$ by
\[N(T) := \big|\{\rho\in \mathcal{Z}: 0<\gamma \le T\}\big|.
\]
Thus in what follows, the phrase \textit{counting with multiplicity} indicates that the underlying set of zeros in each expression is actually the multiset of zeros. On the other hand, since
\[N(T) = \sum_{\substack{\rho \in \mathcal{Z} \\ 0<\gamma \le T}} 1 = \sum_{\substack{\rho \in Z \\ 0<\gamma \le T}}m_\rho ,\]
one can also count zeros in $Z$ \textit{weighted by multiplicity} to recover $N(T)$. For clarity we always indicate which of the sets $Z$ or $\mathcal{Z}$ are being used. We will see that $N^{\circledast}(T)$ defined in \eqref{HorMult} distinguishes the multiplicities and horizontal positions of the zeros by using the property that double sums over zeros both count and weight by multiplicity at the same time. 

\subsection{Sums over zeros} The Riemann-von Mangoldt formula for $N(T)$
 states that

\[
N(T) := \sum_{\substack{\rho\in \mathcal{Z} \\ 0 < \gamma \le T}} 1= \frac{T}{2\pi} \log \frac{T}{2 \pi} - \frac{T}{2\pi} + O(\log T),
\]
from which it immediately follows that

\begin{equation}\label{N(T)2} N(T) = \frac{T}{2\pi}\log T + O(T) \qquad \text{and}\qquad N(T+1)-N(T) \ll \log T. \end{equation}
Using the second estimate in \eqref{N(T)2}, or directly as in \cite[Lemma, Ch. 15]{Dav2000}, we have

\begin{equation}\label{w(u)estimate} \sum_{\rho \in \mathcal{Z}} \frac{1}{1+ (t-\gamma)^2 }\ll \log(|t|+2). \end{equation}

Let $\mathcal{Z}(T)$ denote the multiset of those $\rho \in \mathcal{Z}$ for which $T<\gamma\le 2T$. Then \eqref{N(T)2} also implies

\begin{equation} \label{zerosumnotation}
\sum_{ \rho \in \mathcal{Z}(T)}1 :=\sum_{\substack{ \rho \in \mathcal{Z}\\ T<\gamma\le 2T}}1 = N(2T)-N(T) = \frac{T}{2\pi} \log T +O(T).\end{equation}
Recalling the first zero $\gamma_1 > 14.134\ldots$, we see that \eqref{w(u)estimate} implies 
\begin{equation}\label{w-estimate}
\begin{aligned}
\sum_{\rho, \rho' \in \mathcal{Z}(T)} \frac{1}{1+(\gamma -\gamma')^2} \ll \sum_{\substack{\rho \in \mathcal{Z}\\0<\gamma \le 2T}}\sum_{\rho'\in \mathcal{Z}} \frac{1}{1+(\gamma-\gamma')^2} \ll \sum_{\substack{\rho\in \mathcal{Z}\\0<\gamma \le 2T}}\log\gamma &\le N(2T)\log (2T) \\
&\ll T\log^2T. 
\end{aligned}
\end{equation}

\subsection{Observations concerning horizontal multiplicity}
We may rewrite \eqref{HorMult} as
\[
N^{\circledast}(T) := \sum_{\substack{\rho,\rho' \in \mathcal{Z}(T) \\ \gamma=\gamma'}} 1 = \sum_{\rho \in \mathcal{Z}(T)}H(\gamma),
\quad \text{where} \quad
H(\gamma) := \sum_{\substack{\rho' \in \mathcal{Z}(T)\\ \gamma'=\gamma}} 1.
\]
If we assume RH, then the condition $\gamma =\gamma'$ is equivalent to $\rho=\rho'$, and it follows that $H(\gamma) = m_\rho$, the usual multiplicity of the zero $\rho$. Without RH, $H(\gamma)$ counts with multiplicity the number of zeros on the horizontal line $t=\gamma$. Thus we may make the following observations, which will be helpful in the proof of Theorem \ref{thm3}.
\begin{itemize}
    \item If $H(\gamma)=1$, then on the line $t=\gamma$, there is one simple zero on the critical line.
    \item If $H(\gamma)=2$, then on the line $t=\gamma$, there is either a double zero on the critical line or a pair of two symmetric simple zeros equal distant from the critical line. 
    \item If $H(\gamma)=3$, then on the line $t=\gamma$, there is either a triple zero on the critical line, or a simple zero on the critical line and a pair of two symmetric simple zeros.
\end{itemize}
We remark that if there is a double zero off the critical line with $t=\gamma$, then it is necessary that $H(\gamma)\ge 4$. 

In the proof of Theorem \ref{thm3}, we consider the multiset $\mathcal{S}$ of simple zeros in $\mathcal{Z}(T)$ and the multiset $\mathcal{C}$ of critical zeros in $\mathcal{Z}(T)$. We also consider unions, intersections, and complements of these multisets contained in the universe $\mathcal{Z}(T)$. Note that $\mathcal{Z}(T)$ is partitioned into disjoint multisets of zeros by 
\[ \mathcal{H}(k) := \{\rho \in \mathcal{Z}(T): H(\gamma)=k\}, \qquad k\in \mathbb{N}. \] 
We also see from the classification above for $1\le k \le 3$ that
\[
\mathcal{H}(1) \subset \mathcal{S}\cap\mathcal{C},\quad
\mathcal{H}(2) \subset (\mathcal{S'}\cap\mathcal{C})\cup (\mathcal{S}\cap\mathcal{C'}),\quad
\mathcal{H}(3) \subset (\mathcal{S}\cap \mathcal{C})\cup(\mathcal{S'}\cap \mathcal{C})\cup (\mathcal{S}\cap\mathcal{C'})
= \mathcal{S}\cup\mathcal{C}.
\]
The above relations will be important in the proof. We further notice the following set relations:
\begin{equation}\label{SunionC}
N(\mathcal{S}\cup\mathcal{C}) = N(\mathcal{S}) + N(\mathcal{C}) - N(\mathcal{S}\cap\mathcal{C}) \quad \text{and} \quad N(\mathcal{S}\cup\mathcal{C}) = N(\mathcal{S}\cap\mathcal{C'}) +N(\mathcal{S}\cap\mathcal{C})+N(\mathcal{S'}\cap\mathcal{C}),
\end{equation}
which on combining gives
\begin{equation}\label{SunionC_combined}
N(\mathcal{S}) + N(\mathcal{C}) = 2N(\mathcal{S}\cap\mathcal{C})+N(\mathcal{S'}\cap\mathcal{C})+N(\mathcal{S}\cap\mathcal{C'}).
\end{equation}

Finally, we remark that the condition $\mathbf{C}\ge 1$ in Theorem \ref{thm3} comes from the trivial relation
\[
\sum_{\substack{\rho,\rho' \in \mathcal{Z}(T)\\ \gamma=\gamma'}} 1
~\ge~
\sum_{\substack{\rho,\rho' \in \mathcal{Z}(T)\\ \rho =\rho'}} 1
~=~
\sum_{\rho \in \mathcal{Z}(T)} m_{\rho}
~\ge~
\sum_{\rho \in \mathcal{Z}(T)} 1
~\sim~
\frac{T}{2\pi}\log{T}
\]
as $T\to \infty$, from the asymptotic formula \eqref{zerosumnotation}. 
\subsection{Proof of Theorem \ref{thm3}}

We now assume \eqref{thm3input}.
If $H(\gamma)=1$, using an argument analogous to Montgomery's proof of simple zeros in \cite{Montgomery73}, we have
\[ N(\mathcal{S}\cap\mathcal{C}) \ge \sum_{\substack{\rho \in \mathcal{Z}(T)\\ H(\gamma)=1}} 1 \ge \sum_{\rho \in \mathcal{Z}(T)} (2-H(\gamma))\ge (2-\mathbf{C}+o(1))N(T), \]
which proves \ref{thm3-i} of Theorem \ref{thm3}.
Next, 
\[ \begin{split} (3 -\mathbf{C}+o(1) )N(T)&\le \sum_{\rho \in \mathcal{Z}(T)} (3-H(\gamma))\le 2\sum_{\substack{\rho \in \mathcal{Z}(T)\\ H(\gamma)=1}} 1 +\sum_{\substack{\rho \in \mathcal{Z}(T)\\ H(\gamma)=2}} 1 \\ &
\le 2 N(\mathcal{S}\cap\mathcal{C}) + N(\mathcal{S}'\cap\mathcal{C}) + N(\mathcal{S}\cap\mathcal{C}') = N(\mathcal{S}) + N(\mathcal{C})
\end{split} \]
by \eqref{SunionC_combined}, which proves \ref{thm3-ii} of Theorem \ref{thm3}. For \ref{thm3-iii} of Theorem \ref{thm3}, we have, since $\mathcal{H}(3)$ is disjoint from $\mathcal{H}(1)$ and $\mathcal{H}(2)$ and contained in $\mathcal{S}\cup\mathcal{C}$,
\[ \begin{split} (4-\mathbf{C}+o(1))N(T) &\le \sum_{\rho \in \mathcal{Z}(T)} (4-H(\gamma))\le 3\sum_{\substack{\rho \in \mathcal{Z}(T)\\ H(\gamma) =1 }} 1 +2\sum_{\substack{\rho \in \mathcal{Z}(T)\\ H(\gamma)=2}} 1 + \sum_{\substack{\rho \in \mathcal{Z}(T)\\ H(\gamma)=3}} 1\\ & \le 3N(\mathcal{S}\cap\mathcal{C}) + 2N(\mathcal{S}'\cap\mathcal{C}) + 2N(\mathcal{S}\cap\mathcal{C}') = 2N(\mathcal{S}\cup\mathcal{C}) + N(\mathcal{S}\cap\mathcal{C})
\\& \le
3N(\mathcal{S}\cup\mathcal{C}),
\end{split}\]
by the second equality in \eqref{SunionC}.

\begin{remark} The above discussion is included in Sections 6 and 7 of the expository paper \cite{GS25} written by the second and third authors. In \cite[Section 6]{GS25}, they also give a different proof of Theorem \ref{thm3} which is elementary and completely set theoretical, without using the notion of horizontal multiplicity.
\end{remark}

\begin{remark} For the interested reader, we briefly describe our original approach for proving Theorem~\ref{thm3}, which involved a different initial classification for the solutions of $\gamma = \gamma'$ in $N^{\circledast}(T)$. 
The condition $\gamma = \gamma'$ can occur in three ways: 1) when $\rho = \rho'$, which gives what we call diagonal terms, 2) when $\beta \neq 1/2$, where we get $\rho' = 1-\overline{\rho}$ from two symmetric zeros on a horizontal line on either side of the critical line, which we call symmetric diagonal terms, and 3) when there are 3 or more zeros on a horizontal line from non-symmetric pairs of zeros on that line. Hence we obtain 
\[
N^{\circledast}(T) \ = \
\sum_{\rho \in \mathcal{Z}(T)}m_{\rho}
+
\sum_{\substack{\rho \in \mathcal{Z}(T)\\ \beta \neq \frac{1}{2}}}m_{\rho}
+
\sum_{\substack{\rho,\rho' \in \mathcal{Z}(T)\\ \ \beta+\beta' \neq 1,\, \gamma = \gamma'}} 1,
\]
where the first sum on the right side comes from the diagonal terms, the second is given by the symmetric diagonal terms, and the third is from non-symmetric horizontal terms. Proving Theorem~\ref{thm3} using this initial decomposition is not difficult but is much more tedious than using horizontal multiplicity. For further details, see \cite[Section 6]{GS25}. 
\end{remark}

\section{Unconditional Montgomery Theorem} \label{sec3}

In \cite{Montgomery73}, Montgomery introduced the subject of pair correlation of zeros of the zeta-function. He studied, for $x>0$ and $T\ge 3$,
\[
F(x,T) := \sum_{\substack{\rho, \rho'\in \mathcal{Z} \\ 0<\gamma,\gamma' \le T}} x^{i(\gamma-\gamma')} w(\gamma-\gamma'), \qquad \text{where} \qquad \quad w(u) := \frac{4}{4+u^2}.
\]
While the definition of $F(x,T)$ does not depend on RH, the asymptotic formula Montgomery obtained for $F(x,T)$ when $1\le x\le T$ does depend on RH. This dependence, however, is only because RH allowed the real part of the zeros $\beta$ to be removed out of $F(x,T)$. In \cite{BGST-PC} we retained this dependence to obtain an unconditional form of Montgomery's theorem. To state our result, let
\begin{equation}\label{calF(x,T)}
\mathcal{F}(x,T) :=\sum_{\rho,\rho' \in \mathcal{Z}(T)} x^{\rho-\rho'}W(\rho-\rho'), \qquad \text{where} \qquad W(u) := \frac{4}{4 - u^2}. 
\end{equation}

\begin{MT}\label{MT}
For $x\ge 1$ and $T\ge 3$, we have $\mathcal{F}(x,T)\ge 0$, $\mathcal{F}(x,T)= \mathcal{F}(1/x,T)$, and 
\begin{equation}\label{Mon-1}
\mathcal{F}(x,T) = \frac{T}{2\pi x^2}\left(\log^2T\right) \left(1+ O\left(\frac{1}{\sqrt{\log T}}\right)\right) + \frac{T}{2\pi}\log x + O\left(T\sqrt{\log T}\right),
\end{equation}
uniformly for $1\le x \le T$.
\end{MT}

The statement above has been modified from its original formulation in \cite{BGST-PC}, with two changes. First, note that the sum $\mathcal{F}(x,T)$ above is taken over zeros satisfying $T<\gamma,\gamma'\le 2T$, whereas in \cite{BGST-PC} the sums are taken over zeros satisfying $0<\gamma,\gamma'\le T$. Secondly, and more importantly, the error terms appearing above have been corrected from those in the original statement of the theorem.\footnote[2]{All the applications of Theorem \ref{thm1} in \cite{BGST-PC}, such as Lemma 5 and Lemma 7, remain correct and analogous to those applications obtained using \eqref{Mon-1}. The source of the error came from applying Lemma~8 of \cite{GM87}, which has the same issue but likewise no effect on any consequential results in that paper. We note that Montgomery and Vaughan (to appear) have obtained on RH a more refined version of Theorem \ref{thm1} which deals correctly with the lower order terms which we absorb into a larger error term here.} The result stated in \cite{BGST-PC}, in contrast to \eqref{Mon-1} above, was
\[ \mathcal{F}(x,T) = \Big(\frac{T}{2\pi x^2}\log^2 T + \frac{T}{2\pi} \log x\Big)\left(1 + O\left(\frac{1}{\sqrt{\log T}}\right)\right) . \]
Note that if, for example, $x=c\log T$, then 
\[\mathcal{F}(x,T)=\frac{T}{2\pi}\log\log T + \frac{T}{2\pi} \left( \frac1{c^2} + \log c \right) +O\left(T\frac{\log\log T}{\sqrt{\log T}}\right),\]
and the factor of $c$ on the $T$ term is incorrect. In the corrected version above the error term $O(T\sqrt{\log T})$ which now holds over the whole range $1\le x\le T$ absorbs all these terms in this range. We thank Ram\={u}nas Garunk\v{s}tis and Julija Paliulionyt\.{e} for pointing out the mistake and for suggesting an argument to correct the issue. Below, before we prove \hyperref[MT]{\bf MT} for $\mathcal{F}(x,T)$ when $0< \gamma, \gamma' \le T$, we explain how to correct the proof of \cite[Theorem 1]{BGST-PC}.\\ 

We now correct the proof of \cite[Theorem 1]{BGST-PC} to obtain \hyperref[MT]{\bf MT} for $\mathcal{F}(x,T)$ when $0< \gamma, \gamma' \le T$.
Next, using this result, we prove \hyperref[MT]{\bf MT} for $\mathcal{F}(x,T)$ when $T< \gamma, \gamma' \le 2T$. (This could also be done by modifying each step in the first proof.)

\begin{proof}[Proof of Theorem 1 in \cite{BGST-PC}]

Let $I$ be an interval and
\begin{equation} \label{calF_I} \mathcal{F}_I(x) := \sum_{\substack{\rho,\rho' \in \mathcal{Z} \\ \gamma,\gamma' \in I}} x^{\rho-\rho'}W(\rho-\rho').\end{equation} Thus we see from \eqref{calF(x,T)} that $\mathcal{F}(x,T) = \mathcal{F}_{(T,2T]}(x)$. Then by (2.17) and (2.18) of \cite{BGST-PC} we have, for $1\le x\le T$,
\begin{equation} \label{calF-R}\mathcal{F}_{(0,T]}(x) + O(T^{1/2}) + O(x) = \frac{1}{2\pi}R(x,T), \end{equation}
where by \cite[Lemma 1]{BGST-PC} and following the argument in \cite{Montgomery73}, 
\[ R(x,T) = \int_0^T \left| A_1(x,t) +A_2(x,t) + A_3(x,t) \right|^2\, dt,\]
where 
\[ \begin{split} &A_1(x,t) = \frac{\log(|t|+2)}{x}, \qquad A_2(x,t) =-\sum_{n=1}^\infty \frac{\Lambda(n)}{n^{\frac12+it}}\min\left\{\frac{n}{x} ,\frac{x}{n} \right\}, \\& A_3(x,t) \ll \frac{1}{x} + \frac{x^{\frac12}}{(1+t^2)} + \frac{x^{-\frac52}}{|t|+2}. \end{split}\]
Now, for $1\le x\le T$,
\[ M_1 := \int_0^T\left| A_1(x,t) \right|^2\, dt = \frac{T}{x^2} \left((\log^2 T + O(\log T)\right),\]
and by \cite[Lemma 7]{GM87} or \cite[Lemma B]{Gold81},
\[ M_2 :=\int_0^T\left| A_2(x,t) \right|^2\, dt = T\log x + O\left(T\sqrt{\log T}\right),\]
and
\[ M_3 := \int_0^T\left| A_3(x,t) \right|^2\, dt \ll \frac{T}{x^2} + x.\]
Using Cauchy-Schwarz inequality, and an argument of Ram\={u}nas Garunk\v{s}tis and Julija Paliulionyt\.{e}~\cite{gpemail}, we get
\[ \begin{split} R(x,T) &= \int_0^T \left| A_1(x,t) + A_2(x,t) \right|^2\, dt + O\left(\sqrt{M_1 M_3}\right) + O\left(\sqrt{M_2M_3}\right) + O(M_3)\\&
= M_1 +M_2 + O\left(\left|\int_0^T A_1(x,t) \overline{A_2(x,t)}\, dt\right|\right) + O\left(T\sqrt{\log T}\right) + O\left(\frac{T}{x^2}\log T\right).
\end{split} \]
Since for $n\ge 2$
\[ \int_0^T n^{it}\log (t +2)\, dt \ll \frac{\log T}{\log n},\]
we have 
\[\int_0^T A_1(x,t) \overline{A_2(x,t)}\, dt \ll \frac{\log T}{x} \sum_{n=2}^\infty \frac{\Lambda(n)}{n^{\frac12}\log n}\min\left\{\frac{n}{x} ,\frac{x}{n}\right\} \ll \frac{\log T}{\sqrt{x}\log 2x} \ll \log T.\]
Thus, by \eqref{calF-R}, we have for $1\le x\le T$,
\begin{equation}\label{Fthm1(0,T]}\begin{split} \mathcal{F}_{(0,T]}(x) &= \frac{T}{2\pi x^2}(\log^2 T +O(\log T)) + \frac{T}{2\pi}\log x + O\left(T\sqrt{\log T}\right) \\&=
\frac{T\log^2 T}{2\pi x^2}\left( 1+O\left(\frac1{\sqrt{\log T}}\right)\right) +\frac{T}{2\pi}\log x + O\left(T\sqrt{\log T}\right), \end{split} \end{equation}
where we have degraded the first error term insignificantly to simplify applications to sums over differences of zeros. 
Thus we have obtained \hyperref[MT]{\bf MT} for $\mathcal{F}_{(0,T]}(x)$.
\end{proof}

We now turn to proving \hyperref[MT]{\bf MT} for $\mathcal{F}_{(T,2T]}(x)$.
As $\mathcal{F}_I(x)$ in \eqref{calF_I}, we define $F_I(x)$ for $F(x,T)$ with imaginary part of the zeros in an interval $I$, in similar fashion.
\begin{lem}\label{lem1}
We have, for $x>0$ and $T\ge 3$,
\begin{equation}\label{lem1eq1} \mathcal{F}_{(T,2T]}(x) = \frac{2}{\pi} \int_{-\infty}^{\infty} \left|\sum_{\rho \in \mathcal{Z}(T) } \frac{x^{\rho-1/2}}{1-\left(\rho-(1/2+it)\right)^2} \right|^2 \, dt \end{equation}
and 
\begin{equation}\label{lem1eq2} F_{(T,2T]}(x) = \frac{2}{\pi} \int_{-\infty}^{\infty} \left|\sum_{\rho \in \mathcal{Z}(T)}\frac{x^{i\gamma}}{1+(t-\gamma)^2} \right|^2 \, dt. \end{equation}
\end{lem}
This was proved unconditionally as Lemma 3 in \cite{BGST-PC} with the sum over terms with $0< \gamma \le T$. That proof works line by line for both \eqref{lem1eq1} and \eqref{lem1eq2}. 

\begin{proof}[Proof of \eqref{Mon-1}] From \eqref{lem1eq1} of Lemma \ref{lem1} we see immediately that $\mathcal{F}_{(T,2T]}(x,T)$ is real and non-negative, and since we can interchange $\rho$ and $\rho'$ in \eqref{calF(x,T)}, we have $\mathcal{F}_{(T,2T]}(1/x,T)=\mathcal{F}_{(T,2T]}(x,T)$. Thus, $\mathcal{F}(T^{\alpha},T)$ is real, non-negative, and even. 
From \eqref{calF(x,T)} we have
\[ \begin{split} \mathcal{F}_{(T,2T]}(x,T)
&= \sum_{\rho, \rho' \in \mathcal{Z}(T)} x^{\rho-\rho'}W(\rho-\rho') \\
&= \mathcal{F}_{(0,2T]}(x) - \mathcal{F}_{(0,T]}(x) + O\Bigg(\sum_{\substack{\rho, \rho'\in \mathcal{Z} \\ 0< \gamma\le T\\T< \gamma'\le 2T }} x^{|\beta-\beta'|}|W(\rho-\rho')| \Bigg) \\
&= \frac{T\log^2 T}{2 \pi x^2}\left( 1+O\left(\frac1{\sqrt{\log T}}\right)\right) +\frac{T}{2\pi}\log x \\
&\qquad+ O(T\sqrt{\log T}) + O\Bigg(\sum_{\substack{\rho, \rho'\in \mathcal{Z} \\ 0< \gamma\le T\\T< \gamma'\le 2T }} x^{|\beta-\beta'|}|W(\rho-\rho')|\Bigg),
\end{split}\]
where we have applied \eqref{Fthm1(0,T]} for both $\mathcal{F}_{(0,2T]}(x)$ and $\mathcal{F}_{(0,T]}(x)$. Note that zeros off the critical line in the upper half plane occur in pairs $\rho= \beta +i\gamma$ and $1-\bar{\rho}= 1-\beta +i\gamma $ with $\beta\neq 1/2$. Thus, letting 
\[ \Theta(t) := \max\{\beta : \rho=\beta+ i\gamma, \ 0<\gamma \le t\},\]
we have 
\[\max_{0<\gamma, \gamma'\le T} |\beta -\beta'| = 2(\max_{0<\gamma \le T} \beta) -1=2\Theta(T)-1. \]
Hence
\[ E_1 := \sum_{\substack{\rho, \rho'\in \mathcal{Z} \\ 0< \gamma\le T\\T< \gamma'\le 2T }} x^{|\beta-\beta'|}|W(\rho-\rho')| \ll x^{2\Theta(2T)-1}\sum_{\substack{\rho, \rho'\in \mathcal{Z} \\ 0< \gamma\le T\\T< \gamma'\le 2T }}|W(\rho-\rho')|. \]

Using \eqref{N(T)2} and the bound
\begin{equation}\label{5.5BGSTB24}
|W(\rho -\rho')| \ll w(\gamma -\gamma'),
\end{equation}
which is (5.5) of \cite{BGST-PC}, we see that
\[\begin{split} E_1 &\ll x^{2\Theta(2T)-1} \sum_{\substack{\rho\in \mathcal{Z} \\0<\gamma \le T}} \left(\sum_{\substack{\rho'\in \mathcal{Z} \\T<\gamma' \le 2T}}w(\gamma -\gamma')\right)\\ & \ll x^{2\Theta(2T)-1}\sum_{\substack{\rho\in \mathcal{Z} \\0<\gamma \le T}}\frac{\log T}{(T-\gamma)+1}\\ &\ll x^{2\Theta(2T)-1}\log^3T. \end{split}\]

Korobov \cite{koro, koro2} and Vinogradov \cite{vino} (independently) proved that there is a zero-free region for points $\sigma + it$ in the region in the complex plane determined by $\sigma > 1-\eta(t)\ge \Theta(t)$ with
\[ \eta(t) = \frac{c}{(\log{t})^{2/3}(\log\log{t})^{1/3}}, \qquad \text{for}~ t\ge3,\] 
and some constant $c>0$. (It is well-known that there are no zeros in the region $\sigma \geq 0$ and $0\le t\le 3$.)
Hence $E_1 \ll x^{1-2\eta(2T)}\log^3T$. Thus, for $T^{1/2}\le x \le T$,
\[E_1 \ll x^{1-2\eta(2T)}\log^3T \ll x \exp\left(-c\frac{\log x}{(\log{x})^{2/3}(\log\log{x})^{1/3}}\right)\log^3x \ll x, \] 
while, for $1\le x \le T^{1/2}$,
\[E_1 \ll x^{1-2\eta(2T)}\log^3T \ll T^{1/2}\] since this error term is increasing in $x$ and thus we may take $x=T^{1/2}$ in the previous bound. 
Since for $1\le x \le T$, we have $E_1\ll T^{1/2} + x \ll T\ll T\sqrt{\log T}$,
we conclude
\[ \mathcal{F}_{(T,2T]}(x,T) =\frac{T\log^2 T}{2 \pi x^2}\left( 1+O\left(\frac1{\sqrt{\log T}}\right)\right) +\frac{T}{2\pi}\log x + O\left(T\sqrt{\log T}\right).\]
\end{proof}

\section{Evaluation of a Sum over differences of zeros by Montgomery's Theorem} \label{sec4}

We can use Montgomery's theorem to evaluate many sums over the differences $\rho - \rho'$. If RH is assumed then $\rho - \rho' = i(\gamma-\gamma')$ and we can use the usual real Fourier transform. Thus letting $r$ be in $L^1$, then $\widehat{r}$, the Fourier transform of $r$, is well defined and if $\widehat{r}$ is also in $L^1$ we have
\[ r(\alpha) = \int_{-\infty}^\infty\widehat{r}(t)e(\alpha t )\, dt,\qquad \widehat{r}(t) := \int_{-\infty}^\infty r(\alpha)e(-t\alpha )\, d\alpha, \qquad \text{where}\qquad e(u) := e^{2\pi iu}.
\]

For his results on simple zeros assuming RH, Montgomery \cite{Montgomery73} used the Fej{\'e}r kernel 
\begin{equation}\label{Fejer}
j_F(\alpha) = \max\{1-|\alpha|,0\}, \qquad \widehat j_F(t) = \left(\frac{\sin \pi t}{\pi t}\right)^2. \end{equation}
For improved results, he used the Montgomery-Taylor kernel \cite{Montgomery74,CheerGold} 
\[
j_M(\alpha) = \frac1{1 - \cos{\sqrt{2}}}\left(\frac1{2\sqrt{2}}\sin\left(\sqrt{2}j_F(\alpha)\right) + \frac12j_F(\alpha)\cos\left(\sqrt{2}\alpha\right)\right), \]
with $j_F(\alpha)$ from \eqref{Fejer}, and 
\begin{equation} \label{M-That}
\widehat j_M(w) = \frac1{1 - \cos{\sqrt{2}}}\left(\frac{\sin\left(\frac12(\sqrt{2}-2\pi w)\right)}{\sqrt{2}-2\pi w} + \frac{\sin\left(\frac12(\sqrt{2}+2\pi w)\right)}{\sqrt{2}+2\pi w} \right)^2.
\end{equation}
In the formulas that follow we will take the kernel $j(\alpha)$ to be either $j=j_F$ or $j=j_M$. 

When not assuming RH, we use the complex Fourier transform. Define, for $0\le b \ll 1$,
\begin{equation} \label{hatKb}
\widehat{K_b}(t):= \frac{j(2\pi t)}{\cosh(2\pi b t)}.
\end{equation} 
Then we have, for $z\in \mathbb{C}$,
\begin{equation} \label{TsangK} K_b(z) = \int_{-\infty}^\infty \widehat{K_b}(t) e(zt)\, dt = \int_{-\frac{1}{2\pi}}^{\frac{1}{2\pi}} \frac{j(2\pi t)}{\cosh(2\pi b t)} e^{2\pi i zt}\, dt = \frac1{2\pi}\int_{-1}^{1}\frac{j(\alpha)}{\cosh(b\alpha)}e^{iz\alpha}\, d\alpha.
\end{equation}
Note that $K_b(z)$ is an analytic function for all $z$ since the integral defining it is over a finite interval \cite[2.83]{Titchmarsh1}. Also, note that if $b=0$ then $\widehat{K}_0(t)=j(2\pi t)$. We call $K_b(z)$ a Tsang kernel \cite{Tsang3}.

We now take $z= -i(\rho -\rho')\log T$ in \eqref{TsangK} and have 
\[K_b(-i(\rho -\rho')\log{T}) 
= \frac1{2\pi}\int_{-1}^{1}\frac{j(\alpha)}{\cosh(b\alpha)} T^{\alpha(\rho-\rho')}\, d\alpha .\]
Multiplying by $W(\rho-\rho')$ and summing over pairs of zeros $\rho, \rho'$ with $T<\gamma,\gamma'\le 2T$, we obtain by the definition \eqref{calF(x,T)} that
\begin{equation} \label{K-sum1} \sum_{\rho, \rho' \in \mathcal{Z}(T)}
K_b(-i(\rho -\rho')\log{T}) W(\rho-\rho')
= \frac1{2\pi}\int_{-1}^{1}\frac{j(\alpha)}{\cosh(b\alpha)}\mathcal{F}(T^{\alpha},T) \, d\alpha. \end{equation}

By Montgomery's Theorem \eqref{Mon-1} we see that $\mathcal{F}(T^\alpha,T)$ is an even function of $\alpha$, and the right side of \eqref{K-sum1} is equal to
\[
\left( \frac1{\pi} \int_{0}^{1}\left(\frac{j(\alpha)}{\cosh(b\alpha)}T^{-2\alpha}\log T\left( 1+O\left(\frac1{\sqrt{\log T}}\right)\right) + \frac{\alpha j(\alpha)}{\cosh(b\alpha)}\right)\, d\alpha \right)\frac{T}{2\pi}\log T + O\left(T\sqrt{\log T}\right).
\]
Since
\[\begin{split}&\frac1{\pi} \int_{0}^{1}\frac{j(\alpha)}{\cosh(b\alpha)}T^{-2\alpha}\log T\, d\alpha =\frac1{\pi}\left(\int_{0}^{\frac{\log\log T}{\log T}}+ \int_{\frac{\log\log T}{\log T}}^{1}\right)\frac{j(\alpha)}{\cosh(b\alpha)}T^{-2\alpha}\log T \, d\alpha \\& = \left(j(0) + O\left(\frac{\log\log T}{\log T}\right)\right) \frac1{\pi}\int_{0}^{\frac{\log\log T}{\log T}}T^{-2\alpha}\log T \, d\alpha + O\left( \frac1{\log T}\right) =
 \frac1{2\pi} j(0) + O\left(\frac{\log\log T}{\log T}\right), \end{split}\]
we conclude that
\begin{equation}\label{K-sum2}
\sum_{\rho, \rho' \in \mathcal{Z}(T)} 
K_b(-i(\rho -\rho')\log{T}) W(\rho-\rho')
= \frac1{2\pi}\left(j(0) + 2\int_{0}^{1}\frac{\alpha j(\alpha)}{\cosh(b\alpha)}\, d\alpha + O\left(\frac1{\sqrt{\log T}}\right)\right)\frac{T}{2\pi}\log T.
\end{equation}

\section{The Tsang Kernel} \label{sec5}

To obtain the main properties of the Tsang kernel, we need the following Fourier Transform pair.
\begin{lem}\label{lem2} Let $t$, $x$ be real variables, and take real numbers $y$ and $b$ with $|y|<b$. Then we have the Fourier Transform pair
\begin{equation}\label{Coshtransform} h_{y,b}(t) = \frac{\cosh(2\pi y t)}{\cosh(2\pi b t)} \qquad \text{and} \qquad \widehat{h}_{y,b}(x) = \frac1{b} \left(\frac{\cos\frac{\pi y}{2b}\cosh\frac{\pi x}{2b}}{\cos\frac{\pi y}{b}+\cosh\frac{\pi x}{b}}\right) >0. 
\end{equation}
\end{lem}

\begin{proof}[Proof of Lemma \ref{lem2}] By \cite[Example 3, p. 81]{SteinShakarchi2003}, a residue calculation gives
\[ \int_{-\infty}^\infty \frac{e(-xw)}{\cosh\pi x} \, dx = \frac{1}{\cosh\pi w}, \]
so that $1/\cosh\pi x$ is its own Fourier transform. Letting $x=2bt$ and $w=z/2b$, where $z$ is real, then we have $xw=zt$, and taking the real part of both sides of the above equation, we have
\[ \int_{-\infty}^\infty \frac{\cos(2 \pi zt)}{\cosh (2\pi bt)} \, dt = \frac{1}{2b\cosh\frac{\pi z}{2b}}. \]
While this was derived with $z$ real, if we now take $z=x+iy$ with $|y|<b$ we see the integral is absolutely convergent, and therefore each side of this equation is well defined and equals the same analytic function in the horizontal strip $|y|<b$ in the complex plane. 

Recalling the formulas
\[ \cos z = \cos(x+iy) = \cos x\cosh y -i \sin x \sinh y, \quad \cosh z = \cosh(x+iy) = \cosh x\cos y +i \sinh x \sin y, \]
we have ${\rm Re}(\cos z) = \cos x \cosh y$, ${\rm Re}(\cosh z) = \cosh x \cos y$, and
\[ {\rm Re}\left(\frac{1}{\cosh z}\right) = \frac{{\rm Re}(\cosh\bar{z})}{|\cosh z|^2} = \frac{\cosh x \cos y}{\cos^2y+\sinh^2x}= \frac{\cos y \cosh x}{\frac{1+\cos 2y}{2}+\frac{\cosh 2x -1}{2}} = \frac{2 \cos y \cosh x}{\cos 2y +\cosh 2x}. \] 
Then, by the above formulas and using that $h_{y,b}(t)$ is even in the second line below, we have
\[ \begin{split} \widehat{h}_{y,b}(x) &= \int_{-\infty}^\infty e(-xt)h_{y,b}(t) \, dt \\& = \int_{-\infty}^\infty \frac{\cos(2\pi x t)\cosh(2\pi y t)}{\cosh(2\pi b t)} \, dt \\&
= {\rm Re}\int_{-\infty}^\infty \frac{\cos(2\pi zt)}{\cosh(2\pi b t)} \, dt \\&
= {\rm Re}\frac{1}{2b\cosh\frac{\pi z}{2b}}
\\ & = \frac1{b} \left(\frac{\cos\frac{\pi y}{2b}\cosh\frac{\pi x}{2b}}{\cos\frac{\pi y}{b}+\cosh\frac{\pi x}{b}}\right) >0. \end{split}\]

\end{proof}

\begin{remark} Tsang \cite{Tsang3} uses $h_y(t)$ which is the case $b=1$ of the formula above and obtains $\widehat{h}_y(x)$ from \cite[p. 31, (12)]{Bat54}, but this is for half of the Fourier transform formula above after a change of variable. Thus Tsang's $\widehat{h}_y(x)$ needs to be doubled, but this has no effect on anything in his paper since only the property $\widehat{h}_y(x)>0$ is used there. 
\end{remark}

We next use this lemma to obtain the following lemma, which Tsang \cite[Lemma 1]{Tsang3} proved with $b=1$.

\begin{lem}[K.-M. Tsang]\label{lemTsang} For $0<b\ll 1$, the kernel $K_b(z)$ defined by \eqref{hatKb} and \eqref{TsangK} is an even entire function for $z=x+iy$, $x,y\in \mathbb{R}$, and we have 
\begin{equation}\label{lemTsangPos}
{\rm Re}\,K_b(x+i y)>0 \qquad \text{for all $x$ and $|y|<b$}.
\end{equation}
Also, for $z \in\mathbb{C}$, we have 
\begin{equation}\label{lemTsangBound}
K_b(z) \ll \frac{e^{|{\rm Im}(z)|}}{1+|z|^2}.
\end{equation}
\end{lem}

\begin{proof}[Proof of Lemma \ref{lemTsang}]
By \eqref{TsangK}, we have
\begin{equation} \label{K_b2}K_b(z)= \int_{-\frac{1}{2\pi}}^{\frac{1}{2\pi}} \frac{j(2\pi t)}{\cosh(2\pi b t)} e^{2\pi i zt}\, dt = 2 \int_0^{\frac{1}{2\pi}} \frac{j(2\pi t)}{\cosh(2\pi b t)} \cos(2\pi zt)\, dt.
\end{equation} 
Letting $z=x+iy$ and using ${\rm Re}(\cos(u+iv)) = \cos u \cosh v$ for real $u$ and $v$, we have
\[{\rm Re}\,K_b(x+iy) = 2\int_{0}^{\frac{1}{2\pi}} \frac{j(2\pi t)\cosh(2\pi yt)}{\cosh(2\pi b t)}\cos(2\pi xt)\, dt = 2\int_{0}^{\frac{1}{2\pi}} j(2\pi t)h_{y,b}(t)\cos(2\pi xt)\, dt. \]
Hence, by Lemma \ref{lem2}, for $|y|<b$,
\[
{\rm Re}\, K_b(x+iy) = \bigl(j(2 \pi t)\, h_{y,b}(t)\bigr)^{\widehat{}}\ (x) = \widehat{j(2\pi t)}\ast\widehat{h_{y,b}(t)}(x).
\]
Since $\widehat{j}\ge 0$ by \eqref{Fejer} and \eqref{M-That} and $\widehat{h}_{y,b} >0$ by \eqref{Coshtransform}, we see the convolution is positive and \eqref{lemTsangPos} follows. 
We now prove \eqref{lemTsangBound}. By \eqref{K_b2} with the change of variable $\alpha=2\pi t$, we have
\begin{equation}\label{TsangK2} K_b(z) = \frac{1}{\pi} \int_0^1 J_b(\alpha)\cos(z\alpha)\, d\alpha, \qquad J_b(\alpha) := \frac{j(\alpha)}{\cosh(b\alpha)}.\end{equation}
If $|z|\le 1$, then $ |K_b(z)|\ll \int_0^1 j(\alpha)\, d\alpha \ll 1$,
which proves the estimate in this range.
If $|z|>1$, then by integration by parts 
\[\begin{split} K_b(z) &= \left.\frac{\sin z\alpha}{\pi z}J_b(\alpha) \right|_0^1 - \int_0^1\frac{\sin z\alpha}{\pi z} J_b'(\alpha)\, d\alpha \\ &
= \left.\frac{\cos z\alpha}{\pi z^2} J_b'(\alpha) \right|_0^1 - \int_0^1\frac{\cos z\alpha}{\pi z^2}J_b''(\alpha)\, d\alpha 
\ll \frac{e^{|{\rm Im}(z)|}}{|z|^2}.
\end{split}\]
\end{proof}
Returning to $K_b$ in \eqref{K-sum2},
we see from \eqref{lemTsangPos} of Lemma \ref{lemTsang} that
\begin{equation} \label{Kpos}
{\rm Re}\,K_b(-i(\rho -\rho')\log{T})
= {\rm Re}\,K_b((\gamma-\gamma')\log T -i(\beta-\beta')\log T) > 0 \qquad \text{if} \quad |\beta-\beta'| < \frac{b}{\log T}.
\end{equation}
We want the real part of every term of the sum in \eqref{K-sum2} to be positive (or at least non-negative). This is achieved by assuming in Theorem \ref{thm1} that all the zeros in the sum are in $B_b$, since then $|\beta -\frac12|<\frac{b}{2\log T}$ for all zeros with $T<\gamma \le 2T$, and $\rho\in B_b$ and $\rho'\in B_b$ implies $|\beta-\beta'| < \frac{b}{\log T}$ and \eqref{Kpos} is satisfied.

\begin{lem} \label{lem4} Assume that, for all sufficiently large $T$, all zeros $\rho$ of the Riemann zeta-function with $T<\gamma \le 2T$ are in the region $B_b$ from \eqref{B_0} for some $0<b\ll 1$. Then 
\begin{equation}\label{lem4eq}
 \sum_{\substack{\rho,\rho' \in \mathcal{Z}(T) \\ |\beta -\beta'|< b/\log T}}
\text{\rm Re}\, K_b(-i(\rho -\rho')\log{T}) 
= \frac1{2\pi}\left( j(0) + 2\int_{0}^{1}\frac{\alpha j(\alpha)}{\cosh(b\alpha)}\, d\alpha + O\left( \frac1{\sqrt{\log T}}\right)\right)\frac{T}{2\pi}\log T,
\end{equation}
where every term in the sum above is positive. 
\end{lem}
\begin{proof}[Proof of Lemma \ref{lem4}]
Since $W(\rho-\rho')$ is complex-valued we need to remove this weight when $|\beta-\beta'| < \frac{b}{\log T}$ in \eqref{K-sum2} before applying \eqref{Kpos}. Since $W(0)=1$ and $W(z)-1 = \frac{z^2}{4-z^2}$, we have by \eqref{lemTsangBound},
\[
\begin{aligned}
 \sum_{\substack{\rho,\rho' \in \mathcal{Z}(T)\\ \rho \in B_b,\ \rho'\in B_b}}
&K_b(-i(\rho-\rho')\log{T})
(W(\rho-\rho')-1) \ll
 \sum_{\substack{\rho,\rho' \in \mathcal{Z}(T)\\ |\beta-\beta'| < \frac{b}{\log T}}}\frac{T^{|\beta-\beta'|}}{1+|\rho-\rho'|^2\log^2T} \frac{|\rho-\rho'|^2}{|4-(\rho-\rho')^2|} \\ &
\le \frac{1}{\log^2T} \sum_{\rho,\rho' \in \mathcal{Z}(T)} \frac{T^{b/\log T}}{|4-(\rho-\rho')^2|} \ll \frac{e^b}{\log^2T} \sum_{\rho,\rho' \in \mathcal{Z}(T)}w(\gamma-\gamma')
\ll e^bT \ll T, 
\end{aligned}
\]
where we used \eqref{w-estimate} and \eqref{5.5BGSTB24} in the last line. Hence \eqref{K-sum2} becomes 
\[ \sum_{\substack{\rho,\rho' \in \mathcal{Z}(T) \\ |\beta -\beta'|< b/\log T}}
 K_b(-i(\rho -\rho')\log{T}) 
= \frac1{2\pi}\left( j(0) + 2\int_{0}^{1}\frac{\alpha j(\alpha)}{\cosh(b\alpha)}\, d\alpha + O\left( \frac1{\sqrt{\log T}}\right)\right)\frac{T}{2\pi}\log T,\]
since the assumption $\rho,\rho' \in B_b$ for all zeros in $\mathcal{Z}(T)$ allows us
to insert the weaker condition $|\beta-\beta'| < \frac{b}{\log T}$ into the sum in \eqref{K-sum2} without adding any new terms to the sum, and then this allows us to remove $W(\rho-\rho')$ with an error $\ll T \ll T \sqrt{\log T}$. Taking the real parts proves \eqref{lem4eq}, and by \eqref{Kpos} all the terms in the sum are positive.
\end{proof}

\section{Proof of Theorems \ref{thm1} and \ref{thm2}}
\label{sec6}

\begin{lem} \label{lem5} Assume that, for all sufficiently large $T$, all the zeros of the Riemann zeta-function with $T<\gamma \le 2T$ are in the region $B_b$ for some $0<b\ll 1$. Then, for $j(\alpha) = j_F(\alpha)$ or $j_M(\alpha)$, we have as $T\to \infty$
\begin{equation}\label{Main2} N^{\circledast}(T) \le \left( \mathcal{C}_b(j)+ o(1) \right)\frac{T}{2\pi}\log T, \qquad \text{where} \qquad
 \mathcal{C}_b(j) := \frac{j(0) + 2\int_{0}^{1} \frac{\alpha j(\alpha)}{\cosh(b\alpha)}\, d\alpha }{2\int_{0}^{1}\frac{j(\alpha)}{\cosh(b\alpha)}\, d\alpha} . \end{equation}
\end{lem}
\begin{proof}[Proof of Lemma \ref{lem5}]
Since all the terms in the sum in \eqref{lem4eq} are positive, we have 
\[\sum_{\substack{\rho,\rho' \in \mathcal{Z}(T) \\ \gamma =\gamma'}}
\text{\rm Re}\, K_b(-i(\beta -\beta')\log{T}) \le \sum_{\substack{\rho,\rho' \in \mathcal{Z}(T) \\ |\beta -\beta'|< b/\log T}}
\text{\rm Re}\, K_b(-i(\rho -\rho')\log{T}).
 \]
By \eqref{TsangK2}, we have, since $\cosh(x) \ge 1$ for all real $x$,
\[\text{\rm Re}\, K_b(-i(\beta -\beta')\log{T}) = \frac{1}{\pi} \int_0^1 \frac{j(\alpha)}{\cosh(b\alpha)}\cosh((\beta -\beta')\log{T})\alpha)\, d\alpha \ge \frac{1}{\pi} \int_0^1 \frac{j(\alpha)}{\cosh(b\alpha)}\, d\alpha ,\] and have shown
\[ \left(\frac{1}{\pi} \int_0^1 \frac{j(\alpha)}{\cosh(b\alpha)}\, d\alpha \right)N^{\circledast}(T) \le \sum_{\substack{\rho,\rho' \in \mathcal{Z}(T) \\ |\beta -\beta'|< b/\log T}}
\text{\rm Re}\, K_b(-i(\rho -\rho')\log{T})\]
Hence from Lemma \ref{lem4} we obtain the result in \eqref{Main2}.
\end{proof}
\begin{proof}[Proof of Theorem \ref{thm1}] Letting $b\to 0$ as $T\to \infty$ and using the Fej{\'e}r kernel $j(\alpha) = j_F(\alpha)=1 -\alpha$ for $0\le \alpha \le 1$, we have 
\[\lim_{b\to 0} C_b(j) = \frac{j(0) + 2\int_{0}^{1} \alpha j(\alpha)\, d\alpha }{2\int_{0}^{1}j(\alpha)\, d\alpha}= \frac{1 + 2\int_{0}^{1} \alpha (1-\alpha)\, d\alpha }{2\int_{0}^{1}(1-\alpha)\, d\alpha} = \frac43, \]
and thus 
\[ N^{\circledast}(T) \le (\frac43+o(1))\frac{T}{2\pi}\log T.\]
Theorem \ref{thm1} now follows immediately from Theorem \ref{thm3}.
\end{proof}
\begin{proof}[Proof of Theorem \ref{thm2}] We compute $C_b(j_F)$ and $C_b(j_M)$ using the formulas above and Mathematica to obtain Table \ref{table:2}, as well as the values in Theorem \ref{thm2}. The values given are truncations of the computed values. 
\end{proof}

\begin{center}
\begin{table} [ht!]
\begin{tabular}{||l|| l | l||} 
 \hline
 $b$ & $N^s_0(B_b)|_{j_F}\ $ & $N^s_0(B_b)|_{j_M}$\\ [0.5ex] 
 \hline\hline
 0.001 & .66666& .67250 \\ 
 \hline
 0.2&.66422 & .67019 \\
 \hline
 0.4&.65697 & .66333 \\
 \hline
 0.6&.64509 & .65208 \\
 \hline
 0.8&.62886 & .63670 \\ 
 \hline
 1& .60861 & .61748 \\ 
 \hline
 1.2&.58468 & .59475 \\
 \hline
 1.4&.55743 & .56884 \\
 \hline
 1.6&.52719 & .54003 \\
 \hline
 1.8&.49424 & .50862 \\ 
 \hline
 2&.45887 & .47485 \\ 
 \hline
 2.2&.42130 & .43894 \\
 \hline
 2.4&.38176 & .40109 \\
 \hline
 2.6&.34043 & .36149 \\
 \hline
 2.8&.29747 & .32027 \\ 
 \hline
 3&.25304 & .27760 \\ 
 \hline
 3.2&.20727 & .23357 \\
 \hline
 3.4&.16026 & .18832 \\
 \hline
 3.6&.11214 & .14194 \\
 \hline
 3.8&.06298 & .09451 \\ 
 \hline
4&.01288 & .04612 \\
 \hline
 4.0508&.00022 & .03368 \\
 \hline
 4.187& 0 & .00007 \\
 \hline
\end{tabular}
\vskip .1in
\caption{Comparison of using $j_F$ and $j_{M}$ for $N^s_0(B_b)$.}
\label{table:2}
\end{table}
\end{center}

\section*{Acknowledgments and Funding}
The authors thank the American Institute of Mathematics for its hospitality and for providing a pleasant research environment where the authors met and started this research. S. A. C. B. was supported by NSF DMS-1854398 FRG. A. I. S. was supported by Inamori Research Grant 2024, JSPS KAKENHI Grant Number 22K13895, and Kyushu University International Research Leader Training Program (EBXU0101). C.~T.-B. is partially supported by NSF CAREER DMS-2239681.

\bibliographystyle{alpha}
\bibliography{AHReferences}

\end{document}